\documentclass{article}
\usepackage{amssymb}
\usepackage{epsfig}
\usepackage{amsmath}
\usepackage{amsfonts}

\usepackage{color}
\usepackage{hyperref}

\newcommand{\R}{{{\Bbb R}}}

\newcommand{\N}{{{\Bbb N}}}

\newtheorem{theorem}{\sc Theorem}[section]

\newtheorem{lemma}{\sc Lemma}[section]
\newtheorem{definition}{\sc Definition}[section]
\newtheorem{remark}{\sc Remark}[section]

\newtheorem{example}{\sc Example}[section]

\title{On the discontinuous second--order deviated Dirichlet problem with non--monotone conditions}

\author{Rub\'en Figueroa}

\begin{document}
\maketitle

\begin{center}
 Departamento de
An\'alise Matem\'atica\\
Facultade de Matem\'aticas\\Universidade de Santiago de Compostela, Campus Sur \\
15782 Santiago de
Compostela, Spain\\
  {\bf e-mail:} ruben.figueroa@usc.es
\end{center}

\def\tanh#1{\,{\normalsize tanh}{\,#1}\,}
\def\simbolo#1#2{#1\dotfill{#2}}
\def\qed{\hbox to 0pt{}\hfill$\rlap{$\sqcap$}\sqcup$\medbreak}
\def\theequation{\arabic{section}.\arabic{equation}}
\def\thesection {\arabic{section}}

\begin{abstract}
We provide a new result on the existence of extremal solutions for second--order Dirichlet problems with deviation argument. As a novelty in this work, the nonlinearity need not be continuous or monotone. In order to obtain this new result, we use a generalized monotone method coupled with lower and upper solutions.
\end{abstract}

\noindent
{\bf Primary classification number:} 34B99.

\section{Introduction}

Differential equations with deviated arguments have received a lot of attention in last years. Because of their meaningful interest for modelling real--life processes, see for instance \cite{eng}, \cite{F1}, \cite{har}, \cite{liz} and references therein, nowadays mathematicians work a lot in order to study both qualitative and quantitative properties of these equations.\\

In the present work, we look for a new result on the existence and location of extremal solutions for the following second--order Dirichlet problem with deviated argument:

\begin{equation}\label{p}
\left\{
\begin{array}{ll}
-u''(t)=f(t,u(t),u(\tau(t))), \ t \in I=[0,T]; \\
\\
u(t)=\phi(t), \ t \in [-r,0], \quad u(T)=B,
\end{array}
\right.
\end{equation}
where $\phi$ is a continuous start function and $\tau$ is measurable and such that $\tau(I)\subset [-r,T]$. So this framework includes, in particular, delayed equations. We will search for such solutions inside the set
$$\mathcal{X}=\{u \in \mathcal{C}([-r,T]) \, : \, u_{|I} \in W^{2,1}(I)\}.$$

Problem (\ref{p}) has been studied recently from different points of view: in \cite{fp3} the existence of solutions for a more general problem was obtained via application of Schauder's Theorem, assuming so that nonlinearity $f$ is continuous with respect to all variables. In \cite{fig} continuity was replaced by monotone conditions, and this change allowed us to guarantee the existence of extremal solutions between lower and upper solutions. On the other hand, in \cite{jan2nd} the author studied (\ref{p}) with periodic conditions and with the assumption  $\tau(I) \subset I$. His approach combines continuiuty and one--sided Lipschitz conditions in order to avoid monotonicity, by using a similar approach to that done in \cite{liznieto} for a first--order problem. \\

The main goal of the present paper follows the line of \cite{jan2nd}, but we improve those results in the following ways: we require no continuity conditions, we let deviated argument $\tau$ to take values outside $I$ and one--side Lipschitz constants are replaced by $L^p-$ functions. All these contributions allow us to study a larger class of problems, including classical delay problems, as we will show in last section with some examples. In order to do that, we use a generalized monotone method coupled with lower and upper solutions. \\

This paper is organized as follows: In Section $2$ we provide a result on the existence of a unique solution for problem (\ref{p}) in the case that nonlinearity is Lipschitz--continuous with respect to spatial variables. This result is used later, in Section $3$, for solving auxiliar linear approximations of problem (\ref{p}). Section $2$ is devoted to a comparison result for second--order Dirichlet problems with deviation. In Section $3$ we include the main result on this work. There, we use the work from previous sections in order to define an adequate fixed--point operator in a certain functional interval. This operator will provide us the extremal solutions for (\ref{p}). Finally, in Section $4$ we include some examples of application of our results. \\

In the sequel, $\tau:I \longrightarrow [-r,T]$ is a measurable deviated argument.

\section{Uniqueness result}

As an auxiliar step in the process of linearization of problem (\ref{p}), we provide now a result on the existence of unique solutions for problem (\ref{p}) in the case that nonlinearity $f$ is Lipschitz--continuous.

\begin{theorem}\label{lip} Assume that function $f:I \times \R^2 \longrightarrow \R$ satisfies the following conditions:
\begin{enumerate}
\item[$(H_1)$] For each $x,y \in \R$, the mapping $t \in I \longmapsto f(t,x,y)$ is measurable;
\item[$(H_2)$] For each compact subset $K \subset \R^2$ there exists $\psi_K \in L^1(I,[0,+\infty))$ such that
$$
|f(t,x,y)| \le \psi_K(t) \ \mbox{ for all $(x,y) \in K$;}$$
\item[$(H_3)$]{ There exist nonnegative functions $L_1,L_2$ such that
$$
|f(t,x,y)-f(t,\overline{x},\overline{y})| \le L_1(t) |x-\overline{x}| + L_2(t) |y-\overline{y}|.$$
Moreover, functions $L_1,L_2$ satisfy one of the following:
\begin{enumerate}
\item[$(C_1)$] $L_1,L_2 \in L^{\infty}(I)$ and $||L_1+L_2||_{\infty} < \dfrac{1}{T^2}$;
\item[$(C_2)$] $L_1,L_2 \in L^2(I)$ and $||L_1+L_2||_2 < \left(\dfrac{3}{2T^3}\right)^{1/2}$;
\item[$(C_3)$] $L_1,L_2 \in L^1(I)$ and $||L_1+L_2||_1 < \dfrac{1}{2T}$.
\end{enumerate}
}
\end{enumerate}
Under these conditions, problem (\ref{p}) has a unique solution in $\mathcal{X}$.
\end{theorem}

\noindent {\bf Proof.} First of all, notice that $u \in \mathcal{X}$ is a solution of problem (\ref{p}) if and only if $u$ is a fixed point of the operator $G: \mathcal{C}([-r,T]) \longrightarrow \mathcal{C}([-r,T])$ defined as
$$
\left\{
\begin{array}{ll}
Gu(t)=\phi(t), \quad t \in [-r,0], \\
\\
Gu(t)=\phi(0)+Ct - {\displaystyle \int_0^t (t-s) f(s,u(s),u(\tau(s))) \, ds,} \quad t \in [0,T],
\end{array}
\right.
$$
where
$$
C=\dfrac{1}{T}\left(B-\phi(0)+\int_0^T (T-s) f(s,u(s),u(\tau(s))) \, ds\right).
$$
Notice that operator $G$ is well--defined by virtue of conditions $(H_1)$ and $(H_2)$. So, we will show now that $G$ is a contraction when considering the Banach space $\mathcal{C}([-r,T])$ equipped with its usual supremum norm, $||u||=\max_{t \in [-r,T]} |u(t)|:$ \\

\noindent Given $u,v \in \mathcal{C}([-r,T])$, for each $t \in I$ it is

\begin{align}
|Gu(t)-Gv(t)| & \le \dfrac{t}{T}\left(\int_0^T (T-s) |f(s,u(s),u(\tau(s))) -f(s,v(s),v(\tau(s)))|\, ds\right) \\
&+\int_0^t (t-s) |f(s,u(s),u(\tau(s)))-f(s,v(s),v(\tau(s)))| \, ds,
\end{align}
and then
$$
||Gu-Gv|| \le 2 \left(\int_0^T (T-s) (L_1(s)+L_2(s))\, ds\right)||u-v||.$$
Now notice that:
\begin{enumerate}
\item{If $(C_1)$ holds, then
$$
\int_0^T (T-s) (L_1(s)+L_2(s))\, ds \le ||L_1+L_2||_{\infty} \, \dfrac{T^2}{2}<\dfrac{1}{2},$$
so $G$ is a contraction;}
\item{If $(C_2)$ holds, then
$$
\int_0^T (T-s) (L_1(s)+L_2(s))\, ds \le ||L_1+L_2||_2^2 \, \dfrac{T^3}{3}<\dfrac{1}{2},$$
so $G$ is a contraction;}
\item{If $(C_3)$ holds, then
$$
\int_0^T (T-s) (L_1(s)+L_2(s))\, ds \le T \, ||L_1+L_2||_1 <\dfrac{1}{2},$$
so $G$ is a contraction;}
\end{enumerate}
By application of Banach's fixed--point theorem, operator $G$ has a unique fixed point which is the unique solution of problem (\ref{p}) in $\mathcal{X}$. \qed

\section{Maximum principle and main result}

The first step on developing our generalized monotone method is to obtain a comparison result for second--order Dirichlet problems with delay. The proof of this result follows the line of \cite[Theorem 2.2]{jiangwei} and \cite[Theorem 2.1]{liznieto}.

\begin{lemma}\label{maximo} Assume that $p \in \mathcal{X}$ and there exist nonnegative functions $L_1,L_2$ satisfying the following:
\begin{enumerate}
\item[$(M_1)$] $-p''(t)+L_1(t)p(t)+L_2(t)p(\tau(t)) \ge 0$ for a.a. $t \in [0,T]$;
\item[$(M_2)$] $p(0)\ge 0$, $p(T) \ge 0$;
\item[$(M_3)$] $0 \le p(t) \le p(0)$ for all $t \in [-r,0]$;
\item[$(M_4)$]{ functions $L_1,L_2$ satisfy one of the following conditions:
\begin{enumerate}
\item[$\hat{(C_1)}$] $L_1,L_2 \in L^{\infty}(I)$ and $||L_1+L_2||_{\infty} < \dfrac{2}{T^2};$
\item[$\hat{(C_2)}$] $L_1,L_2 \in L^{2}(I)$ and $||L_1+L_2||_2 < \dfrac{\sqrt{2}}{T};$
\item[$\hat{(C_3)}$] $L_1,L_2 \in L^{1}(I)$ and $||L_1+L_2||_1 < \dfrac{1}{T}\;$
\end{enumerate}
}
\end{enumerate}
Then $p(t) \ge 0$ for all $t \in [-r,T]$.
\end{lemma}

\noindent {\bf Proof. } Assume that there exists $t_0$ in $(0,T)$ such that $p(t_0) <0$. We consider two cases: \\
Case $I$: $p(t) \le 0$ for all $t \in [0,T]$ and $p$ is not identically $0$ in that interval. In this case, $p(0)=p(T)=0=p(t)$ for $t \in [-r,0]$, and condition $(M_1)$ provides that
$$
p''(t) \le L_1(t)p(t)+L_2(t)p(\tau(t)) \le 0 \ \mbox{ for a.a. $t \in (0,T)$},$$
so $p \equiv 0$, a contradiction. \\
Case $II$: there exist $t_1,t_2 \in (0,T)$ such that $p(t_1) >0$ and $p(t_2)<0$. In this case, there exists $t_3 \in (0,T)$, $t_4 \in [0,T]$ satisfying
$$
p(t_3)=\min_{t \in [-r,T]} p(t) <0, \quad p'(t_3)=0,$$
$$
p(t_4)=\max_{t \in [-r,T]} p(t) >0.$$
So,
\begin{equation}\label{desig1}
p''(t) \le (L_1(t)+L_2(t)) p(t_4), \ \mbox{ for a.a. $t \in [0,T]$}.\end{equation}
If $t_4 < t_3$ then we take $s \in (t_4,t_3)$ and integrate (\ref{desig1}) from $s$ to $t_3$ to obtain
$$
-p'(s) \le p(t_4) \int_s^{t_3} (L_1(t) + L_2(t)) \, dt,$$
and now integrating from $t_4$ to $t_3$:
$$
p(t_4) \le p(t_4)-p(t_3) \le p(t_4) \int_{t_4}^{t_3}\int_s^{t_3} (L_1(t)+L_2(t)) \, dt \, ds,$$
and then
\begin{equation}\label{int}
1 \le \int_{t_4}^{t_3}\int_s^{t_3} (L_1(t)+L_2(t)) \, dt \, ds.\end{equation}
Now notice that we can bound integral in (\ref{int}) by the following numbers: \\
If condition $\hat{(C_1)}$ holds, then
$$
\int_{t_4}^{t_3}\int_s^{t_3} (L_1(t)+L_2(t)) \, dt \, ds \le \dfrac{T^2}{2} \, ||L_1+L_2||_{\infty} <1,$$
which is a contradiction. \\
If condition $\hat{(C_2)}$ holds, then
$$
\int_{t_4}^{t_3}\int_s^{t_3} (L_1(t)+L_2(t)) \, dt \, ds \le \dfrac{T^2}{2} \, ||L_1+L_2||_2^2 <1,$$
again a contradiction. \\
If condition $\hat{(C_3)}$ holds, then
$$
\int_{t_4}^{t_3}\int_s^{t_3} (L_1(t)+L_2(t)) \, dt \, ds \le  T\, ||L_1+L_2||_1 <1.$$
On the other hand, if $t_3 < t_4$, then we reason like above, now integrating from $t_3$ to $s$ $(t_3<s<t_4)$ and then from $t_3$ to $t_4$ to get the same contradictions. \qed

\begin{remark} Notice that for $\dfrac{3}{4} \le T$ each condition $(C_i)$ implies $\hat{(C_i)}$, $i=1,2,3$, and for
$0<T \le \dfrac{3}{4}$ we have $(C_1) \Rightarrow \hat{(C_1)},$ $(C_3) \Rightarrow \hat{(C_3)}$ and $\hat{(C_2)} \Rightarrow (C_2)$.
\end{remark}

Now we introduce a Lemma on extremal fixed--points for monotone operators in ordered sets. This result will play an essential role in our future argumentations.

\begin{lemma} \cite[Theorem 1.2.2]{heikkila} Let $Y$ be a subset of an ordered metric space $X$, $[a,b]$ a nonempty interval in $Y$ and $G:[a,b] \longrightarrow [a,b]$ a nondecreasing operator. If $\{Gx_n\}_{n=1}^{\infty}$ converges in $Y$ whenever $\{x_n\}_{n=1}^{\infty}$ is a monotone sequence in $[a,b]$, then operator $G$ has the greatest, $x^*$, and the least, $x_*$, fixed point in $[a,b]$. Moreover, we have that

$$
x_*=\min \{x \, : \, Gx \le x\}, \quad x^*=\max \{x \, : \, x \le Gx\}.$$

\end{lemma}

We define now what we mean by lower and upper solutions for problem (\ref{p}).

\begin{definition}\label{lu} We say that $\alpha, \beta \in \mathcal{X}$ are, respectively, a lower and an upper solution for problem (\ref{p}) if the compositions
$$
t \in [0,T] \longmapsto f(t,\alpha(t),\alpha(\tau(t))), \quad t \in [0,T] \longmapsto f(t,\beta(t),\beta(\tau(t)))$$
are measurable and the following conditions hold:
$$
\left\{
\begin{array}{ll}
-\alpha''(t) \le f(t,\alpha(t),\alpha(\tau(t))) \ \mbox{ for a.a. $t \in [0,T]$}, \\
\\
\alpha(t) \le \phi(t), \ \phi(t)-\alpha(t) \le \phi(0)-\alpha(0) \mbox{ for all $t \in [-r,0]$}, \quad \alpha(T) \le B;
\end{array}
\right.
$$
$$
\left\{
\begin{array}{ll}
-\beta''(t) \ge f(t,\beta(t),\beta(\tau(t))) \ \mbox{ for a.a. $t \in [0,T]$}, \\
\\
\beta(t) \ge \phi(t), \ \beta(t)-\phi(t) \le \beta(0)-\phi(0) \ \mbox{ for all $t \in [-r,0]$}, \quad \beta(T) \ge B.
\end{array}
\right.
$$
\end{definition}

\begin{theorem}\label{main} Assume that there exist $\alpha,\beta \in \mathcal{X}$ which are, respectively, a lower and an upper solution for problem (\ref{p}), with $\alpha(t) \le \beta(t)$ for all $t \in [-r,T]$, and put $$[\alpha,\beta]=\{\gamma \in \mathcal{C}([-r,T]) \, : \, \alpha(t) \le \gamma(t) \le \beta(t) \ \mbox{ for all $t \in [-r,T]$}\}.$$ Assume, moreover, that the following conditions hold:
\begin{enumerate}
\item[$(H_1)$] For each $\gamma \in [\alpha,\beta]$ the composition $t \in I \longmapsto f(t,\gamma(t),\gamma(\tau(t)))$ is measurable;
\item[$(H_2)$] There exists $\psi \in L^1(I,[0,\infty))$ such that for a.a. $t \in I$, all $x \in [\alpha(t),\beta(t)]$ and all $y \in [\alpha(\tau(t)),\beta(\tau(t))]$ we have $|f(t,x,y)| \le \psi(t)$;
\item[$(H_3)$] There exist nonnegative functions $L_1,L_2$ such that for a.a. $t \in I$
$$
f(t,\overline{x},\overline{y})-f(t,x,y) \ge -L_1(t) (\overline{x}-x) - L_2(t) (\overline{y}-y)$$
whenever $\alpha(t) \le x\le \overline{x} \le \beta(t)$, $\alpha(\tau(t)) \le y \le \overline{y} \le \beta(\tau(t)).$ \\

Moreover, if $T \ge \dfrac{3}{4}$ then functions $L_1,L_2$ satisfy one of the conditions $(C_1)$, $(C_2)$, $(C_3)$ and if $0<T<\dfrac{3}{4}$ then they satisfy one of the following: $(C_1)$, $\hat{(C_2)}$, $(C_3).$
\end{enumerate}
In these conditions, problem (\ref{p}) has the extremal solutions in $[\alpha,\beta]$.
\end{theorem}

\noindent {\bf Proof.} We define an operator $G:[\alpha,\beta] \longrightarrow [\alpha,\beta]$ as follows: for each $\gamma \in [\alpha,\beta],$ $G\gamma$ is the unique solution of the problem
\begin{equation}\label{aux}
\left\{
\begin{array}{ll}
-u''(t)+L_1(t)(u(t)-\gamma(t)) +L_2(t)(u(\tau(t))-\gamma(\tau(t))) =f(t,\gamma(t),\gamma(\tau(t))), \ \mbox{ for a.a. $t \in I$},\\
\\
u(t)=\phi(t) \ \mbox{ for all $t \in [-r,0]$}, \quad u(T)=B.
\end{array}
\right.
\end{equation}

\noindent {\it Step $1$: Operator $G$ is well--defined from $[\alpha,\beta]$ to $\mathcal{X}$.} We can rewrite the differential equation in (\ref{aux}) in the form $-u''(t)=\widetilde{f}(t,u(t),u(\tau(t)),$ where
$$
\widetilde{f}(t,x,y)=-L_1(t)x - L_2(t) y + f(t,\gamma(t),\gamma(\tau(t)))+L_1(t)\gamma(t)+L_2(t)\gamma(\tau(t)),$$
so problem (\ref{aux}) satisfies conditions of Theorem \ref{lip} and then it has a unique solution.
\\

\noindent {\it Step $2$: Operator $G$ is nondecreasing.} Let $\gamma_1,\gamma_2 \in [\alpha,\beta]$ such that $\gamma_1(t) \le \gamma_2(t)$ for all $t \in [-r,T]$. We will show now that $G\gamma_1(t) \le G\gamma_2(t)$ for all $t$.\\
First of all, notice that $G\gamma_1,G\gamma_2$ satisfy for a.a. $t \in I$ that
$$
-G\gamma_1''(t)+L_1(t)(G\gamma_1(t)-\gamma_1(t))+L_2(t)(G\gamma_1(\tau(t))-\gamma_1(\tau(t))))=f(t,\gamma_1(t),\gamma_1(\tau(t))),$$
$$
-G\gamma_2''(t)+L_1(t)(G\gamma_2(t)-\gamma_2(t))+L_2(t)(G\gamma_2(\tau(t))-\gamma_2(\tau(t))))=f(t,\gamma_2(t),\gamma_2(\tau(t))),$$
and, moreover,
$$
G\gamma_1(t)=G\gamma_2(t)=\phi(t) \ \mbox{ for all $t \in [-r,0]$,} \quad G\gamma_1(T)=G\gamma_2(T)=B.$$
So, by virtue of condition $(H_3)$ we obtain that function $p=G\gamma_2-G\gamma_1$ satisfies:
$$
\left\{
\begin{array}{ll}
-p''(t)+L_1(t)p(t)+L_2(t)p(\tau(t)) \ge 0, \ \mbox{ for a.a. $t \in I$}, \\
\\
p(0)=p(T)=0=p(t) \ \mbox{ for all $t \in [-r,0]$}.
\end{array}
\right.
$$
Then, by application of Lemma \ref{maximo} we obtain that $p \ge 0$ on $[-r,T]$, so $G\gamma_1(t) \le G\gamma_2(t)$ for all $t \in [-r,T]$, which implies that operator $G$ is nondecreasing. \\

\noindent {\it Step $3$: $G([\alpha,\beta]) \subset [\alpha,\beta]$.} Because of being $\alpha$ a lower solution for problem (\ref{p}), we have that $p=G\alpha-\alpha$ satisfies:
$$
\left\{
\begin{array}{ll}
-p''(t)+L_1(t)p(t)+L_2(t)p(\tau(t))\ge 0 \ \mbox{ for a.a. $t \in I$}, \\
\\
p(t)=\phi(t)-\alpha(t)\ge 0, \quad p(t)=\phi(t)-\alpha(t)\le \phi(0)-\alpha(0)=p(0) \ \mbox{ for all $t \in [-r,0]$}, \\
\\
p(T)=\phi(T)-\alpha(t)=T-\alpha(t) \ge 0.
\end{array}
\right.
$$
So, by application of Lemma (\ref{maximo}) we conclude that $G\alpha(t) \ge \alpha(t)$ for all $t \in [-r,T]$. In the same way we prove that $G\beta \le \beta$. This and monotonicity of operator $G$ imply that $G([\alpha,\beta]) \subset [\alpha,\beta]$. \\

\noindent {\it Step $4$: $\{G\gamma_n\}_{n=1}^{\infty}$ converges in $\mathcal{C}([-r,T])$ whenever $\{\gamma_n\}_{n=1}^{\infty}$ is a monotone sequence in $[\alpha,\beta]$.} Put $z_n=G\gamma_n \in \mathcal{C}([-r,T])$. Because of the monotonicities of sequence $\{\gamma_n\}_{n=1}^{\infty}$ and operator $G$ we obtain that $\{z_n\}_{n=1}^{\infty}$ is a monotone sequence in $[\alpha,\beta]$, so it has it pointwise limit, say $z$.\\
As $\{z_n(t)\}_{n=1}^{\infty}$ is a constant sequence for $t \in [-r,0]$, $\{z_n\}_{n=1}^{\infty}$ converges to $z$ uniformly in that interval. Now, fixed $t \in I$, by virtue of Mean Value Theorem there exists $c=c(n) \in [0,T]$ such that
$$
z_n'(t)=\dfrac{B-\phi(0)}{T} + \int_c^t \widetilde{f}(s,z_n(s),z_n(\tau(s))) \, ds,$$
so
$$
||z_n'||_{\infty} \le \dfrac{|B-\phi(0)|}{T} + ||\psi||_1 + ||L_1+L_2||_1 ||\beta-\alpha||_{\infty}.$$
Then, the sequence $\{z_n'\}_{n=1}^{\infty}$ is uniformly bounded on $I$, so $\{z_n\}_{n=1}^{\infty}$ converges to $z$ in $\mathcal{C}(I)$. \\

\noindent {\it Step $5$: Problem (\ref{p}) has the extremal solutions in $[\alpha,\beta]$.} By application of Lemma \ref{maximo}, operator $G$ has in $[\alpha,\beta]$ the extremal fixed points, $u^*$, $u_*$. We will show now that $u^*$, $u_*$ correspond, respectively, with the greatest and the least solution of problem (\ref{p}) in $[\alpha,\beta]$. First of all, it is clear that each fixed point of $G$ is also a solution of (\ref{p}). On the other hand, if $u$ is a solution of (\ref{p}) between $\alpha$ and $\beta$ then $u$ also solves (\ref{aux}), and the uniqueness of solution of problem (\ref{aux}) provides that $Gu=u$. Then, $u_* \le u \le u^*$, so $u_*$, $u^*$ are the extremal solutions of (\ref{p}) in $[\alpha,\beta]$. \qed

\begin{remark} Among all conditions in the previous result, perhaps $(H_1)$ is the most difficult to check in practise when working with discontinuous nonlinearities. The reader is referred to \cite{appell}, \cite{cp2}, \cite{hasrzy} for a larger disquisition about this. In those references, some results for guaranteing $(H_1)$ are provided. On the other hand, and roughly speaking, condition $(H_3)$ forces $f$ not to have downwards discontinuities. In Section $4$ we introduce a Lemma which will be useful on checking this condition in examples.
\end{remark}

\section{Examples of application}

In this section we show some examples of application of Theorem \ref{main}. As far as we know, the following problems can be studied with no result in the literature. As a technical support in order to check condition $(H_3)$, we begin by introducing a simple lemma.

\begin{lemma}\label{crec} Let $\{x_k\}_{k=1}^{\infty}$ a strictly increasing sequence of real numbers and assume that $f:\R \longrightarrow \R$ is such that $f_{|(x_{k-1},x_{k})} \in \mathcal{C}^1(x_{k-1},x_k)$ for all $k=1,2,\ldots$ Assume moreover that the following conditions hold for each $k \in \N$:

\begin{enumerate}

\item[$(i)$]
$$
\lim_{x \to x_k^-}f(x) \le f(x_k) \le \lim_{x \to x_k^+} f(x);$$

\item[$(ii)$] there exist $M_k=\inf \{f'(x) \, : \, x \in (x_{k-1},x_k)\}$ and $M=\inf \{M_k \, : \, k \in \N\}.$

\end{enumerate}

In these conditions, the function $g:x \in \R \longmapsto g(x)=f(x) + |M| x$ is nondecreasing.
\end{lemma}

\noindent {\bf Proof.} Fixed $k \in \N$ we have that $g$ is differentiable in $(x_{k-1},x_k)$ and $g'(x) \ge M_k + |M| \ge 0$ for all $x \in (x_{k-1}, x_{k})$, so $g$ is nondecreasing in that interval. On the other hand, condition $(i)$ provides that $f$ is nondecreasing at point $x_k$, and so $g$. \qed

\begin{example} Consider the following problem with delay ($[\cdot]$ means integer part):
\begin{equation}\label{ex2}
\left\{
\begin{array}{ll}
-u''(t)=[t \, u(t)] - \dfrac{1}{9} u(t-1) \, \sin{\left(\dfrac{u(t-1)\pi}{2[|u(t-1)|]+2}\right)}\equiv f(t,u(t),u(t-1)) \ \mbox{ a.e. on $I=[0,2]$,} \\
\\
u(t)=\phi(t)=\cos{\dfrac{\pi}{2}\, t}, \ t \in [-1,0], \quad u(2)=\dfrac{\pi}{4}.
\end{array}
\right.
\end{equation}

Notice that function $f$ is discontinuous with respect to its three variables and it is non--monotone with respect to its functional one. We will show that problem (\ref{ex2}) has extremal solutions between adequate lower and upper solutions. \\

Consider the functions
\begin{equation}\label{alpha}
\alpha(t)=0 \ \mbox{ for all $t \in [-1,2]$,}
\end{equation}

\begin{equation}\label{beta}
\beta(t)=\left\{
\begin{array}{lr}
\cos{\dfrac{\pi}{2}\, t}, &\mbox{ if $t \in [-1,0]$,} \\
\\
1-t(t-2), &\mbox{ if $t \in [0,2]$,}
\end{array}
\right.
\end{equation}

We will prove that $\alpha$ and $\beta$ are, respectively, a lower and an upper solution for problem (\ref{ex2}). First of all, notice that $\alpha(t) \le \phi(t) \le \beta(t)$ for all $t \in [-1,0]$ and that the compositions $t \in [0,2] \longmapsto f(t, \alpha(t),\alpha(t-1))$ and $t \in [0,2] \longmapsto f(t,\beta(t),\beta(t-1))$ have at most a null set of discontinuity points, so they are both measurable. Now, taking into account for $t \in [0,2]$ it is $t-t^2(t-2) <3$, we have for a.a. $t \in [0,2]$ that
$$
f(t,\alpha(t),\alpha(t-1)) = f(t,0,0)= 0=-\alpha''(t),$$
$$
f(t,\beta(t),\beta(t-1)) = [t-t^2(t-2)]- \dfrac{1}{9} \beta(t-1) \sin{\left(\dfrac{\beta(t-1)\pi}{2[|\beta(t-1)|]+2}\right)} \le 2 = -\beta''(t).$$

Then, $\alpha$ and $\beta$ are lower and upper solutions for our problem, which moreover satisfy $\alpha \le \beta$ on $[-1,2]$. \\

Reasoning as above, if $\gamma \in [\alpha,\beta]$, with $[\alpha,\beta]$ defined as in Theorem \ref{main}, then the composition $t \in [0,2] \longmapsto f(t,\gamma(t),\gamma(t-1))$ is measurable, so condition $(H_1)$ in that Theorem holds. \\

On the other hand, for $t \in [-1,2]$ we have that $0 \le \alpha(t) \le \beta(t) \le 3$, so condition $(H_2)$ holds with $\psi \equiv 4$. \\

We are going now to check $(H_3)$. First, notice that for $t \in [0,2]$ function $f(t,\cdot,y)=[tx]$ is nondecreasing, so we can take $L_1 \equiv 0$. On the other hand, for $t \in [0,2]$ and $x \in [\alpha(t),\beta(t)]$ we have that $f(t,x,\cdot)$ is discontinuous exactly at integer numbers, with
$$
\lim_{y \to k^-} f(t,x,y)=-\dfrac{1}{9} y \sin\left(\dfrac{\pi}{2}\right)$$
and
$$
f(t,x,k)=\lim_{y \to k^+} f(t,x,y)=-\dfrac{1}{9} y \sin\left(\dfrac{k}{k+1} \, \dfrac{\pi}{2}\right),$$
for $k \in \mathbb{Z}$. So, condition $(i)$ in Lemma \ref{crec} is satisfied for $y \ge 0$ (in particular, inside the interval $[\alpha(t-1),\beta(t-1)]$). Moreover, for $y \in (k-1,k)$, $k \ge 1$, it is
$${\displaystyle
\dfrac{\partial}{\partial y} f(t,x,y)=\dfrac{-1}{9} \left(1+\dfrac{y \pi}{2k}\right) \cos{\left(\dfrac{y \pi}{2k}\right)} \ge -\dfrac{1}{9} \left(1+\dfrac{\pi}{2}\right),}$$
and by application of Lemma \ref{crec} condition $(H_3)$ is satisfied with $L_2 \equiv \dfrac{1}{9} \left(1+\dfrac{\pi}{2}\right).$ Notice that $\dfrac{1}{9} \left(1+\dfrac{\pi}{2}\right) \approx 0.2857,$ so nor condition $(C_1)$ nor $(C_3)$ are satisfied. Ho\-we\-ver, $\sqrt{2} \, L_2 < \dfrac{\sqrt{3}}{4},$ so condition $(C_2)$ holds. \\

As all the conditions in Theorem \ref{main} hold, we conclude that problem (\ref{ex2}) has the extremal solutions between $\alpha$ and $\beta$.

\end{example}

\begin{example} Fix $0 < k < \dfrac{1}{10}$, let $\varphi:\R \longrightarrow \R$ the function defined by
$$
\varphi(0)=0, \quad \varphi(x)=\dfrac{k}{n}-kx \ \mbox{ if $|x| \in \left(\dfrac{1}{n+1},\dfrac{1}{n}\right] \bigcup (n,n+1]$, $n=1,2,\ldots$},$$
and consider the following problem with both delay and advance:
\begin{equation}\label{ex1}
\left\{
\begin{array}{ll}
-u''(t)=\sin(t)+\varphi(u(\sqrt{1-t}))+\dfrac{1}{5\sqrt{t}}u(\sqrt{t}) \equiv f(t,u(\sqrt{1-t}),u(\sqrt{t})) \ \mbox{ for a.a. $t \in I=[0,1]$}, \\
\\
u(0)=u(1)=0.
\end{array}
\right.
\end{equation}
Notice that in problem (\ref{ex1}) the nonlinearity $f$ is discontinuous and non--monotone with respect to its second variable in a countable set. Moreover, it blows up when $t \to 0$. \\

Consider now the functions $\alpha(t)=t^2-t=-\beta(t)$ for $t \in [0,1]$. We will show that $\alpha$ and $\beta$ are, respectively, a lower and an upper solution for problem (\ref{ex1}). \\

First, notice that if $\gamma$ is a continuous function then $\varphi \circ \gamma$ is discontinuous at most in a countable set, so $\varphi \circ \gamma$ is measurable. On the other hand, $\alpha(0)=\alpha(1)=0=\beta(0)=\beta(1)$, and for a.a. $t \in I$ we have that
$$
f(t,\alpha(\sqrt{1-t}),\alpha(\sqrt{t}))=\sin(t)+\varphi(1-t-\sqrt{1-t})+\dfrac{\sqrt{t}-1}{5} \ge -1-1/5 \ge-2=-\alpha''(t),$$
$$
f(t,\beta(\sqrt{1-t}),\beta(\sqrt{t}))=\sin(t)+\varphi(\sqrt{1-t}-(1-t))+\dfrac{1-\sqrt{t}}{5} \le 1 + k/2 + 1/5 \le 2=-\beta''(t).$$
Then, $\alpha,\beta$ are, repectively, a lower and an upper solution for (\ref{ex1}). \\

Now notice that $\dfrac{-1}{4} \le \alpha(t) \le \beta(t) \le \dfrac{1}{4}$ for all $t \in I$ and for $x,y \in \left[\dfrac{-1}{4},\dfrac{1}{4}\right]$ we have
$$
|f(t,x,y)| \le \psi(t)= \sin(t) + \dfrac{k}{2} + \dfrac{1}{20\sqrt{t}},$$
with $\psi \in L^1(I)$, so condition $(H_2)$ holds. \\

Finally, condition $(H_3)$ is satisfied with $L_1(t)=k$ and $L_2(t)=\dfrac{1}{5\sqrt{t}}$, $t \in (0,1]$, and we have that
$$
||L_1+L_2||_1=k+\dfrac{2}{5} < \dfrac{1}{2} \ \mbox{ if $k < \dfrac{1}{10}$.}$$

By application of Theorem \ref{main}, problem (\ref{ex1}) has the extremal solutions in $[\alpha,\beta]$.
\end{example}

\begin{remark} Notice that, in problem \ref{ex1}, there exists no constant $\hat{L}$ such that $$y \in [\alpha({\sqrt{t}}),\beta({\sqrt{t}})] \longmapsto f(t,x,y)+\hat{L}y$$ is nondecrasing, so our improvement (from constant one--sided Lipschitz conditions to $L^p-$ ones) become essential in this case.
\end{remark}

\bibliographystyle{amsplain} 

\end{document}